\documentclass{article}
\usepackage{color,amsmath,amssymb,latexsym,amsfonts}
\usepackage{mathrsfs}

\renewcommand{\abstractname}{Abstract.}
\renewcommand\abstract{\hfil\break\topsep=0pt\partopsep=0pt\parsep=0pt\itemsep=0pt\relax
\trivlist\item[\hskip\labelsep
{\bfseries\abstractname}]\if!\abstractname!\hskip-\labelsep\fi}

\newcommand{\email}[1]{{(e-mail: #1)}}

\def\keywordname{{\bfseries Key words:}}
\def\keywords#1{\par\addvspace\baselineskip\noindent\keywordname\enspace
\ignorespaces#1}

\def\subclassname{{\bfseries Mathematics Subject Classification (1991):} }
\def\subclass#1{\par\addvspace\medskipamount\noindent\subclassname\
\ignorespaces#1}

\def\title#1{\hfil\break\hfil\break
\hfil\break\par\addvspace\baselineskip\noindent
\ignorespaces{\LARGE\bf#1}\hfil\break}

\def\author#1{\par\addvspace\baselineskip\noindent
\ignorespaces{\large\bf#1}}

\def\institute#1{\par\addvspace\baselineskip\noindent
\ignorespaces{\small#1}\hfil\break}

\newtheorem{theorem}{Theorem}

\newtheorem{lemma}[theorem]{Lemma}

\newtheorem{proposition}[theorem]{Proposition}

\def\bds{\begin{displaystyle}}
\def\eds{\end{displaystyle}}

\newcommand{\F}{{\mathbb F}}

\begin{document}

\title{On nonexistence of non-constant volatility in the Black-Scholes formula\footnote{This
work was supported by the Australian Research Council.\\
Date: September 2004}}

\author{Hamza, K. and Klebaner, F.C.}
\institute{School of Mathematical Sciences, Monash University}
\email{kais.hamza@sci.monash.edu.au,
fima.klebaner@sci.monash.edu.au}

\begin{abstract}
We prove  that if the Black-Scholes formula holds with the spot
volatility for call options with all strikes, then the volatility
parameter is constant. The proof relies some result
on semimartingales (Theorem \ref{specmartg}) of independent interest.
\end{abstract}



\subclass{60G44, 60H30, 90A09}

\section{Introduction}
One of the most important applications of stochastic dynamics has
occurred in finance as a model for evolution of prices of stocks
and their options. This model in its simple form is described by a
randomly perturbed  exponential growth. If $S_t$ denotes the price
of stock at time $t$, then its evolution is given by the
stochastic differential equation
\begin{equation} dS_t=\mu S_t dt +\sigma S_t dB_t, \label{BSeqn}
\end{equation}  with $B_t$ denoting the Brownian motion process. The strength of
the random perturbation is determined by the parameter $\sigma$,
which is known as the volatility of the stock in finance. The
above model was used by Merton, Black and Scholes to find the
price of an option on stock, such as an agreement to buy the stock
at some future time $T$ for the specified at time $t<T$ price $K$.
Their formula states that the price of such an option at time $t$
is given by
\begin{equation}
C_t=S_t\Phi\Big(\frac{\log
\frac{S_t}{K}+(r+\frac{\sigma^2}{2})(T-t))}{\sigma\sqrt{T-t}}\Big)-
Ke^{-r(T-t)}\Phi\Big(\frac{\log \frac{S_t}{K}+(r-\frac{\sigma^2}{2}
)(T-t))}{\sigma\sqrt{T-t}}\Big), \label{BSform}
\end{equation}  where $\Phi$ denotes the standard
normal distribution function. Remarkably, the parameter $\mu$ does
not enter the formula, but $\sigma$ does, as well as $r$, the
riskless rate available in a  savings account. The Black-Scholes
formula is widely used in financial markets and risk management.

 The volatility of the stock is seen to be  the
parameter in the quadratic variation of the return on the stock
process, $\sigma^2=d[R,R]_t/dt$, where $dR_t=dS_t/S_t$, see
(\ref{BSeqn}). It is widely believed and experimentally verified
that stocks do not have a constant volatility, rather this
parameter varies with time, see e.g. \cite{Spok}, \cite{Gold},
\cite{Stein}. In spite of this fact, the Black-Scholes options
pricing formula (\ref{BSform}) is still used with some adjustments
to the volatility. Thus the question of existence of a model with
non-constant volatility in which the Black-Scholes formula remains
valid is of interest in financial mathematics as well as in
practical applications. A non-constant volatility model has the
form
\begin{equation}
dS_t=\mu S_t dt +\theta_t S_t dB_t, \label{BSnonconst}
\end{equation}  where $\theta_t$ is a
function of time that can be  random,  for example, it can be a
function of stock $S_t$, as well as include other independent
source of randomness. There is a large literature on non-constant
volatility models, both deterministic and stochastic (see e.g.
Fouque et.al. (2001)). In the next paragraph we need to recall the
basic facts on options pricing (see e.g. Shiryaev (1999)). We can
assume in without loss of generality that the riskless interest
rate $r=0$, otherwise work with discounted prices $S_te^{-rt}$.

The First Fundamental theorem of asset pricing states that a model
does not admit arbitrage if and only if there exists an equivalent
probability measure $Q$ such that $S_t$ is a $Q$-martingale.
 The price at time $t$ of a call option
that pays $(S_T-K)^+$ at time $T$ is given by
\begin{equation}
C_t=E_Q((S_T-K)^+|{\mathcal F}^S_t), \label{option}
\end{equation}
where $E_Q$ is the expectation  under $Q$ and ${\mathcal F}^S_t$
is the $\sigma$-field generated by the process $S_u, u\le t$.

Next we comment on   recent work to reconcile the non-constant
volatility with the use of the Black-Scholes formula  by
constructing the so-called stochastic implied volatility (SIV)
market models. A market model is a model that returns
Black-Scholes option prices.  Such models are important in
practical applications for calculations of non-standard options.
The idea in SIV models is to use equation (\ref{BSnonconst})
together with
  the volatility surface   $\sigma(t,T,K)$. These $\sigma(t,T,K)$ are implied by
  the Black-Scholes option prices by equating the theoretical
  price in (\ref{option})
$C_t(T,K)=E_Q((S_T-K)^+|{\mathcal F}^S_t)$ with the Black-Scholes
price
  (\ref{BSform}) with $\sigma=\sigma(t,T,K)$.
 The spot volatility $\theta_t$ is obtained from the
  volatility surface $\sigma(t,T,K)$ by $\theta_t=\sigma(t,t,S_t)$.
 This approach results in a
nonlinear system of stochastic differential equations with delay,
see Brace et.al. (2001), also Schonbucher (1998), Carr (2000),
Brace et.al.  (2002). While existence of  SIV is still an open
problem, the purpose of this note is to show that, unfortunately,
non-constant volatility models are not compatible with the
Black-Scholes formula.

Comment that the ``Black-Scholes" equation with non-constant
volatility $dS_t = \theta_tS_tdB_t$ holds for a wide class of
positive martingales.
 If $S_t$ is a positive martingale with
$P(S_t>0)=1$ for any $t\le T$ and  the predictable representation
property holds then there exists a process $\theta_t$, such that
$dS_t=\theta_t S_t dB_t$.
  We prove the statement even for wider class of processes and do not assume any
dynamics on $S_t$, namely
    if for a positive   semimartingale $S_t$ and for all values of $K$ and three
different maturities the Black-Scholes formula holds  with some
adapted $\theta_t$, then $\theta_t$ must be constant. This
$\theta_t$ can be taken as any adapted functional of the spot
volatility, including the spot volatility itself.

\section{Results}
Let $(\Omega, {\mathcal F}, \F , Q)$ be the filtered probability
space, with the general conditions,  supporting a Brownian motion
$B_t$. Let for a constant $\sigma>0$, $$dZ_t = \sigma
Z_tdB_t,\quad Z_0=z_0$$ and
\begin{equation}
C(T,t,K,\sigma,z) = \mathbb{E}[(Z_T-K)^+|Z_t=z].\label{main}
\end{equation}
\begin{theorem}\label{black}
Let $S_t$ and $\theta_t$ be two adapted processes such that
$\theta_0 = \sigma$ and $S_0=z_0$. Assume that $S_t$ is strictly
positive, that ${\mathcal{F}}_0$ is trivial and that there exist
three equidistant terminal times, $T_1<T_2<T_3$ such that, for all
$K$ and all $t\leq T_i$, $i=1,2,3$
$$\bds \mathbf{E}[(S_{T_i}-K)^+|{\mathcal{F}}_t] =
C(T_i,t,K,\theta_t,S_t)\eds.$$ Then $\theta^2_t = \sigma^2$ for
all $t\le T_1$.
\end{theorem}

The proof relies on the following result.

\begin{theorem}\label{specmartg}
Let $M_t$ and $X_t$ be two  semimartingales, $M_t$ is strictly
positive. Assume that $M_t$, $M_tX_t$ and $M_tX_t^2$ are local
martingales.
Then $X_t\equiv X_0$.
\end{theorem}

Note that Theorem \ref{black} can be formulated for three terminal
non-equidistant times $T_i$, $i=1,2,3$, and in this case Theorem
\ref{specmartg} can be generalized to the case when $M_t$,
$M_tX_t$ and $M_tX_t^{\alpha}$, $\alpha>1$ are local martingales.
For  clarity of presentation we use equidistant times.

\bigskip
\noindent {\bf Proof} of Theorem \ref{specmartg}. Denote by $X^c$
the continuous martingale component of $X$, by $M^c$ that of $M$,
and by $V$ any optional compensator of $X$. By the It\^o formula
we find that
\begin{equation*}
\begin{aligned}
M_tX_t&=M_0X_0+\int_0^tM_{s-}dX_s+\int_0^tX_{s-}dM_s+\langle
X^c,M^c\rangle_t+ \sum_{s\le t}\triangle X_s\triangle M_s.
\\
(MX)_tX_t&=(MX)_0X_0+\int_0^tX_{s-}d(XM)_s+\int_0^tM_{s-}X_{s-}dX_s
\\
&\qquad
+\int_0^tM_{s-}d\langle X^c,X^c\rangle_s+\int_0^tX_{s-}d\langle M^c,X^c\rangle_s
\\
&\qquad +\sum_{s\le t}\triangle(MX)_s\triangle X_s.
\end{aligned}
\end{equation*}
Extracting martingales and using the   assumptions of the theorem
it follows that the random processes
\begin{equation*}
\begin{aligned}
L^{(1)}_t&=\int_0^tM_{s-}dV_s+\langle X^c,M^c\rangle_t+ \sum_{s\le
t}\triangle X_s\triangle M_s,
\\
L^{(2)}_t&=\int_0^tM_{s-}X_{s-}dV_s
\\
&\qquad
+\int_0^tM_{s-}d\langle X^c,X^c\rangle_s+\int_0^tX_{s-}d\langle X^c,M^c\rangle_s
\\
&\qquad +\sum_{s\le t}\triangle(MX)_s\triangle X_s
\end{aligned}
\end{equation*}
are local martingales. Therefore the process
\begin{eqnarray}
L^{(2)}_t-\int_0^tX_{s-}dL^{(1)}_s&=&\int_0^tM_{s-}d\langle
X^c,X^c\rangle_s\nonumber
\\
&&\qquad +\sum_{s\le t}\left[\triangle(MX)_s\triangle
X_s-X_{s-}\triangle M_s\triangle X_s\right]\nonumber
\\
& =&\int_0^tM_{s-}d\langle X^c,X^c\rangle_s+\sum_{s\le
t}M_s(\triangle X_s)^2 \label{L2L1}
\end{eqnarray}
is a local martingale. Since $M_t>0$, it is an increasing process.
Therefore it is zero. It follows that $\triangle X_t=0$ and
$\int_0^tM_{s-}d\langle X^c,X^c\rangle_s=0$.

Next we repeatedly use the following simple argument. If for a
positive function $f$ and an increasing function $g$,
$\bds\int_0^tf(s)dg(s)\equiv0\eds$ then $\bds
g(t)-g(0)=\int_0^tf(s)^{-1}dh(s)\equiv0\eds$ where $\bds
h(t)=\int_0^tf(s)dg(s)\eds$. In other words $g$ is constant.

Applying this argument to $\int_0^tM_sd\langle
X^c,X^c\rangle_s=\int_0^tM_{s-}d\langle X^c,X^c\rangle_s$, we see
that $\langle X^c,X^c\rangle_t\equiv0$, that is $X$ is a
continuous process of finite variation. Since
$$M_tX_t =M_0X_0+\int_0^t X_sdM_s +\int_0^t M_{s-}dX_s,$$
it follows that $\bds\int_0^t M_sdX_s=\int_0^t M_{s-}dX_s\eds$ is
a continuous local martingale, and since it is of finite variation
$\bds\int_0^t M_sdX_s\equiv 0\eds$. Applying the above argument to
$\bds\int_0^t M_sdY_s\eds$, where $Y$ is the variation process of
$X$, we complete the proof. \hfill$\Box$

\bigskip

The proof  of Theorem \ref{black} is broken into a number of
propositions. All of them assume the conditions and notations of
Theorem \ref{black}.

\begin{proposition}\label{Prop3}
 $S_t$
is a martingale and $\mathbf{E}[S_t^2]<+\infty$ for all $t$.
\end{proposition}
\noindent {\bf Proof} Take any $i$. By condition (\ref{main})
$S_{T_i}$ is integrable and
\begin{eqnarray*}
 \mathbf{E}[S_{T_i}|{\mathcal{F}}_t] &=&
\mathbf{E}[(S_{T_i}-0)^+|{\mathcal{F}}_t] =
C(T_i,t,0,\theta_t,S_t) \\
&=&\mathbf{E}[(Z_{T_i}-0)^+|Z_t=z]_{\sigma=\theta_t,z=S_t}\\
&=&\mathbf{E}[Z_{T_i}|Z_t=z]_{\sigma=\theta_t,z=S_t}\\
&=&\left.Z_t\right|_{\sigma=\theta_t,z=S_t} = S_t,
\end{eqnarray*}
 which proves that $S_t$ is a martingale.

Also $\bds\mathbf{E}[(S_{T_i}-K)^+] =
\mathbf{E}[(S_{T_i}-K)^+|{\mathcal{F}}_0] =
C(T_i,0,K,\theta_0,S_0) = C(T_i,0,K,\sigma,z_0) =
\mathbf{E}[(Z_{T_i}-K)^+]\eds$. Using Lemma \ref{lemma1}, we
deduce that $S_{T_i}$ and $Z_{T_i}$ have the same distribution.
The integrability of $Z_{T_i}^2$ induces that $S_{T_i}^2$ and by
Jensen's inequality that of $S_t^2$.\hfill$\Box$

\begin{proposition}
$\bds\left(S_t^2e^{\theta_t^2(T_i-t)}\right)_{t\leq T_i}\eds$ is a
martingale, for  each $i=1,2,3$.
\end{proposition}
\noindent {\bf Proof} The proof is based on  the representation of
$\mathbf{E}[X^2|{\mathcal{G}}]$ in terms of integrals
$\int_0^{+\infty}\mathbf{E}[(X-K)^+|{\mathcal{G}}]dK$ and
$\int_0^{+\infty}\mathbf{E}[(X+K)^-|{\mathcal{G}}]dK$ given in
Lemma \ref{lemma2}. Using  this lemma we obtain
\begin{eqnarray*}
\lefteqn{\mathbf{E}[S_{T_i}^2|{\mathcal{F}}_t]}\\ & = &
2\int_0^{+\infty}\mathbf{E}[(S_{T_i}-K)^+|{\mathcal{F}}_t]dK\\
&& +
2\int_{-\infty}^0\left(\mathbf{E}[(S_{T_i}-K)^+|{\mathcal{F}}_t]-\mathbf{E}[S_{T_i}|{\mathcal{F}}_t]+K\right)dK\\
& = & 2\left[\int_0^{+\infty}C({T_i},t,K,\sigma,z)dK +
\int_{-\infty}^0\left(C({T_i},t,K,\sigma,z)-z+K\right)dK\right]_{\sigma=\theta_t,z=S_t}\\
& = &
2\left[\int_0^{+\infty}\mathbf{E}[(Z_{T_i}-K)^+|Z_t=z]dK\right.\\
&& + \left.\int_{-\infty}^0\left(\mathbf{E}[(Z_{T_i}-K)^+|Z_t=z]-\mathbf{E}[Z_{T_i}|Z_t=z]+K\right)dK\right]_{\sigma=\theta_t,z=S_t}\\
& = &
\left[\mathbf{E}[Z_{T_i}^2|Z_t=z]\right]_{\sigma=\theta_t,z=S_t}\\
& = & \left[z^2e^{\sigma^2(T_i-t)}\right]_{\sigma=\theta_t,z=S_t}\\
& = & S_t^2e^{\theta_t^2(T_i-t)}
\end{eqnarray*}
Thus for each $i$,
$\bds\left(S_t^2e^{\theta_t^2(T_i-t)}\right)_{t\leq T_i}\eds$ is a
martingale of the Doob's form.\hfill$\Box$

\begin{proposition}
$\bds\left(\theta_t^2\right)_{t\leq T_1}\eds$ is a semimartingale.
\end{proposition}
\noindent {\bf Proof} Since $S>0$, $\bds e^{\theta_t^2(T_2-T_1)} =
\frac{S_t^2e^{\theta_t^2(T_2-t)}}
{S_t^2e^{\theta_t^2(T_1-t)}}\eds$ is a semimartingale in
$t\in[0,T_1]$. Thus $\bds\left(\theta_t^2\right)_{t\leq T_1}\eds$
is also a semimartingale. \hfill$\Box$

\bigskip

\noindent {\bf Proof of Theorem \ref{black}}

 Let $\bds M_t =
S_t^2e^{(T_1-t)\theta_t^2}\eds$, $X_t = e^{(T_2-T_1)\theta_t^2} =
e^{(T_3-T_2)\theta_t^2}$. It is easy to see using the previous
propositions that $M_t$ and $X_t$ satisfy the conditions of
Theorem \ref{specmartg}. By its conclusion $X_t \equiv X_0$ and
$\theta^2_t = \sigma^2$ for all $t\le T_1$, which completes the
proof. \hfill$\Box$

\bigskip
\noindent {\bf Remarks}

1. Note that although in Theorem \ref{black} we do not assume any
dynamics on $S_t$, the assumptions of theorem imply that $S_t$ is
a strictly positive martingale, see Proposition \ref{Prop3}. If in
addition, the predictable representation property with respect to
$B_t$ holds, then $S_t$ can be represented as a stochastic
volatility model $dS_t=h_t S_t dB_t$, see e.g. \cite{Klebaner},
p.286.

2. Let  $dS_t=h_t S_t dB_t$. While we may primarily think of
$\theta_t$ in Theorem \ref{black} as the volatility
$\theta_t=h_t$, Theorem \ref{black} applies   for  any adapted
functional  of $h_u$ and $S_u$, $u\leq t$, for example, the
average volatility on $[0,t]$ is given by   $ \theta_t =
\sqrt{\frac{1}{t}\int_0^th^2_udu}  $.

3.   Let $dS_t=h_t S_t dB_t$, where $h_t$ is a deterministic
function. Then
$S_t=S_0\exp(\int_0^th_sdB_s-\frac{1}{2}\int_0^th^2_sds)$ has a
lognormal distribution and the price of the option is given by the
Black-Scholes formula $  \mathbf{E}[(S_{T}-K)^+|{\mathcal{F}}_t] =
C(T,t,K,\theta(t,T),S_t)$ with
 $ \theta^2(t,T) =
\frac1{T-t}\int_t^Th_u^2du $, see e.g. [9]. This example does not
contradict Theorem \ref{black} since in this case $\theta_t$
depends also on $T$.

\section{Auxiliary Results}
\begin{lemma}\label{lemma1}
If $\mathbf{E}[|X|]<\infty$, then $G(K) =
\mathbf{E}[(X-K)^+|{\mathcal{G}}]$ is a convex function
($G^\prime_-$ and $G^\prime_+$ are increasing, respectively left
and right-continuous and $\{K:G^\prime_-(K)\neq G^\prime_+(K)\}$
is countable) with
$$G^\prime_-(K) = -\mathbf{P}[X\geq K|{\mathcal{G}}]\;\mbox{ and }\;G^\prime_+(K) = -\mathbf{P}[X>
K|{\mathcal{G}}].$$
\end{lemma}
\noindent {\bf Proof} The proof easily follows from the fact that,
for $\varepsilon>0$,
$$(x-(K+\varepsilon))^+-(x-K)^+ = \left\{\begin{array}{ll}
0&x\leq K\\
-(x-K)&K<x<K+\varepsilon\\
-\varepsilon&x\geq K+\varepsilon
\end{array}\right.$$
and
$$(x-(K-\varepsilon))^+-(x-K)^+ = \left\{\begin{array}{ll}
0&x\leq K-\varepsilon\\
x-K+\varepsilon&K-\varepsilon<x<K\\
\varepsilon&x\geq K
\end{array}\right.$$

Note that $\bds1_{K<X<K+\varepsilon} = 1_{K-\varepsilon<X<K} =
0\eds$ for $\varepsilon$ small enough.

%
\hfill$\Box$

%
%

\begin{lemma}\label{lemma2}
If $\mathbf{E}[X^2]<+\infty$, then
\begin{eqnarray*}
\lefteqn{\frac12\mathbf{E}[X^2|{\mathcal{G}}]}\\
& = & \int_0^{+\infty}\mathbf{E}[(X-K)^+|{\mathcal{G}}]dK +
\int_0^{+\infty}\mathbf{E}[(X+K)^-|{\mathcal{G}}]dK\\
& = & \int_0^{+\infty}\mathbf{E}[(X-K)^+|{\mathcal{G}}]dK +
\int_{-\infty}^0\left(\mathbf{E}[(X-K)^+|{\mathcal{G}}]-\mathbf{E}[X|{\mathcal{G}}]+K\right)dK
\end{eqnarray*}
In particular, if $X$ is non-negative,
$$\mathbf{E}[X^2|{\mathcal{G}}] =
2\int_0^\infty\mathbf{E}[(X-K)^+|{\mathcal{G}}]dK.$$
\end{lemma}
\noindent {\bf Proof} It is easily checked that for any $x$,
$$x^2 = 2\int_0^{+\infty}(x-K)^+dK + 2\int_0^{+\infty}(x+K)^-dK.$$
The result immediately follows. \hfill$\Box$

\end{document}